\documentclass[11pt]{amsart}
\usepackage{geometry}                
\usepackage{graphicx}
\usepackage{amssymb}
\usepackage{epstopdf}
\title{The upper bound of the Mertens function from the viewpoint of statistical mechanics}
\author{Rong Qiang Wei}
\address{College of Earth and Planetary Sciences, University of Chinese Academy of Sciences, Beijing, PRC, 100049}
\email{wrq1973@ucas.edu.cn}
\date{}

\begin{document}
\maketitle

\begin{abstract}

   We provide some upper bounds for the Mertens function ($M(n)$: the cumulative sum of the M$\ddot{\mathrm{o}}$bius function) by an approach of statistical mechanics, in which the M$\ddot{\mathrm{o}}$bius function is taken as a particular state of a modified one-dimensional (1D) Ising model without the exchange interaction between the spins. Further, based on the assumptions and conclusions of the statistical mechanics, we discuss the problem that $M(n)$ can be equivalent to the sum of an independent random sequence. It holds in the sense of equivalent probability, from which another two upper bounds for the $M(n)$ can be inferred. Besides, if $M(n)$ is a measured quantity, its upper bound is $\sqrt{\frac{B}{\alpha} n}$ ($B$ is constant) with a probability $>1-\alpha$ ($0<\alpha<1$) from the view point of the energy fluctuations in the canonical ensemble.
      
\end{abstract}

{\hspace{2.2em}\small Keywords:}

{\hspace{2.2em}\tiny Mertens function, Upper bound, Statistical mechanics, Ising model, Energy fluctuation}

\section{Introduction}

The Mertens function $M(n)$ is defined as the cumulative sum of the M$\ddot{\mathrm{o}}$bius function $\mu(k)$ for all positive integers $n$,

\begin{equation}\label{eq1}
M(n)=\sum_{k=1}^{n}\mu(k)
\end{equation}
where the M$\ddot{\mathrm{o}}$bius function $\mu(k)$ is defined as follows for a positive integer $k$ by

\begin{equation}\label{eq2}
\mu(k)  = \left\{ {\begin{array}{*{20}{c}}
1\\
0\\
{{{( - 1)}^m}}
\end{array}\begin{array}{*{20}{c}}
{k = 1{\mbox{\hspace{14em} }}}\\
{{\mbox{if \ }} k {\mbox{\ is divisible by a prime square\hspace{2em}}}}\\
{{\mbox{if\  }} k {\mbox{\  is the product of }}m{\mbox{ distinct primes}}}
\end{array}} \right.
\end{equation}

In most cases, $M(n)$ can be extended to real numbers as follows:

\begin{equation}
M(x)= \sum_{n\le x}\mu(n)
\end{equation}

The upper bound of $M(n)$ is always of concern to people, and many results have been presented. The following are just a few examples. The famous one is an old conjecture, "Mertens conjecture" , proposed that $\vert M (n) \vert < n^{1/ 2}$ for all $n$, which was disproved by Odlyzko and te Riele (1985).  Wei (2016) showed that Mertens conjecture is not true in a statistical point of view. Based on an assumption that $\mu (n)$ is an independent random sequence, Wei (2016) discussed the upper bound of the $M(n)$, and found that the inequality (\ref{eq4}) for $M(n)$ holds with a probability of $1 - \alpha$, 

\begin{equation}\label{eq4}
M(n) \le  \sqrt{\frac{6}{\pi^2}} K_{_{\frac{\alpha}{2}}}\sqrt{n}
\end{equation}
where,
$$
\int\limits_{ - {K_{\alpha /2}}}^{{K_{\alpha /2}}} {\frac{1}{{\sqrt {2\pi } }}} \exp ( - \frac{{{t^2}}}{2}){\rm{d}}t = 1 - \alpha 
$$

or, the following inequality (\ref{eq5}) holds with a probability $p>1-\alpha$ 

\begin{equation}\label{eq5}
M(n) \le  \frac{\sqrt {6/{\pi ^2}}}{\sqrt{\alpha}} \sqrt{n}
\end{equation}

Without taking $\mu (n)$ as an independent random sequence, Wei (2017) conjectured that the upper bound of the $M(n)$ should have a similar formula to (\ref{eq4}) or (\ref{eq5}) from three facts. 

Besides the examples above, MacLeod (1967; 1969) showed that

\begin{equation}
\left | M(x)\right | \le \frac{x+1}{80}+\frac{11}{2} \hspace{5em} (x\ge 1)
\end{equation}
 
El Marraki (1995) proved that,

\begin{equation}
\left | M(x) \right | \le \frac{0.002969}{(\log x)^{1/2}}x \hspace{5em} (x\ge 142194)
\end{equation}

and

\begin{equation}
\left | M(x) \right | \le \frac{0.6437752}{\log x}x \hspace{5em} (x > 1)
\end{equation}

Ramar$\acute{e}$ (2013) showed that,

\begin{equation}
\left | M(x) \right | \le \frac{0.0146\log x-0.1098}{(\log x)^2} x \hspace{5em} (x \ge 464402)
\end{equation}

In this paper,  the upper bound of the $M(n)$ will be estimated based on a modified 1D Ising model in the statistical mechanics. Although it is possible that the results are not as precise as those mentioned above, we can study the upper bound of $M(n)$ and the similar problems from a new perspective.

\section{1D Ising Model and its partition function}\label{sec2}

The 1D Ising model in the statistical mechanics is composed of a chain of $n$ spins, in which each spin interacts only with its two nearest neighbors and with an external magnetic field $h$. This nearest-neighbor Ising model is defined in terms of the following Hamiltonian (total energy) (eg., Huang, 1987),

\begin{equation}\label{eq6}
\mathcal{H} =  -\frac{1}{2}J\sum\limits_{\left\langle {i,j} \right\rangle} s_i s_j-\xi h\sum_{i} s_i
\end{equation}
where ${\left\langle {i,j} \right\rangle}$ means that the sum is carried out over all the nearest-neighbor pair of spins $(i,j)$. $s_i=+1$, $s_i=-1$  are the two possible states. $J$ denotes the exchange interaction energy. $\xi $ is the magnetic moment of each spin. 

For our purpose, we modify this Ising model in which $s_i=0$ besides $s_i=\pm1$ according to the value of $\mu(k)$.
When the periodic boundary condition ($s_0=s_{_n}$) is imposed, the partition function for such a modified Ising model is,

\begin{equation}\label{eq7}
Q_{_n}=\sum_{s_i=+1,0,-1}\exp(-\frac{\mathcal{H}}{kT})=\sum_{s_i=+1,0,-1}\exp\left [\frac{1}{kT}\left (\frac{1}{2}J\sum\limits_{\left\langle {i,j} \right\rangle} s_i s_j+\xi h\sum_{i} s_i\right )\right ]
\end{equation}
where $\sum\limits_{\left\langle {i,j} \right\rangle} s_i s_j=s_0s_1+s_1s_2+s_2s_3+\ldots$+$s_{_{n-1}}s_{{_n}}$, and $\sum_{i} s_i=s_1+s_2+s_3+\ldots$+$s_{{_n}}$. $k$ is the Boltzmann constant and $T$  is the temperature. $\sum_{s_i=+1,0,-1}$ is to be understood to extend over all possible states of the model system, namely,
 
\begin{equation}\label{eq8}
\sum_{s_i=+1,0,-1} \thicksim \sum_{s_1=+1,0,-1}  \sum_{s_2=+1,0,-1}\sum_{s_3=+1,0,-1}\ldots \sum_{s_n=+1,0,-1}
\end{equation} 

To evaluate $Q_{_n}$, we define a matrix $\bf{P}$ as follows referring to Huang (1987),

\begin{equation}\label{eq11}
 {\bf{P}} = \left( {\begin{array}{*{20}{c}}
1&1&1\\
y&{xy}&{\frac{y}{x}}\\
{\frac{1}{y}}&{\frac{1}{{xy}}}&{\frac{x}{y}}
\end{array}} \right)
\end{equation}
where $x=\exp (\frac{J}{2kT})$,$y=\exp(\frac{\xi h}{kT})$.

Similar to Huang (1987),  it can also be proven the following eq. (\ref{eq12}), 

\begin{equation}\label{eq12}
Q_{_n}=\mbox{tr\ } \bf{P}^n=\lambda_1^n+\lambda_2^n+\lambda_3^n
\end{equation}
where $\lambda_1$,$\lambda_2$, and $\lambda_3$ are the three eigenvalues of $\bf{P}$, and they are three roots of the equation ${\rm{det}}\vert \bf{P}-\lambda\vert =0$.

In our case as follows,  the exchange interaction energy is not taken into account, i.e., $J=0$, and $x=1$. Then we have 

\begin{equation}\label{eq11+1}
 {\bf{P}} = \left( {\begin{array}{*{20}{c}}
1&1&1\\
y&{y}&{y}\\
{\frac{1}{y}}&{\frac{1}{{y}}}&{\frac{1}{y}}
\end{array}} \right)
\end{equation}

From eq. (\ref{eq11+1}), we can obtain $\lambda_1=\lambda_2=0$, and $\lambda_3=(1+y+\frac{1}{y})$ from ${\rm{det}}\vert \bf{P}-\lambda\vert =0$.  Finally,    

\begin{equation}\label{eqQn}
Q_{_n}=(1+y+\frac{1}{y})^n=\left[1+\exp(\frac{\xi h}{kT})+\exp(-\frac{\xi h}{kT})\right]^n
\end{equation}

\section{The upper bound of the $M(n)$}\label{result}

Here we investigate a particular state of the 1D Ising model system in the section \ref{sec2}, namely the arrangement of spins is: $s_0=s_{_n}=\mu(n)$, $s_1=\mu(1)=+1$, $s_2=\mu(2)=-1$, $s_3=\mu(3)=-1$, $s_4=\mu(4)=0$, $\ldots$, $s_{_{n-1}}=\mu(n-1)$.  In the case of $x=1$, the Hamiltonian $\mathcal{H'}$ for this arrangement is,

 \begin{equation}\label{eq13}
\mathcal{H'} = -\xi h\sum_{i} s_i
\end{equation}

According to the statistical mechanics, the probability $p$ for this special state is,

\begin{equation}\label{eq14}
p=\frac{\exp(-\frac{\mathcal{H'}}{kT})}{Q_{_n}}=\frac{\exp(\frac{\xi h}{kT}\sum_{i} s_i)}{Q_{_n}}=\frac{\exp(\beta\sum_{i} s_i)}{Q_{_n}}=\frac{\exp\left [\beta M(n)\right ]}{Q_{_n}}
\end{equation}

Since $p\leq 1$, we can obtain an upper bound for $M(n)$ as the follows,

\begin{equation}\label{eq15}
M(n)\leq \frac{1}{\beta}\log Q_{_n}=\left\{\frac{1}{\beta}\log \left[1+\exp(\beta)+\exp(-\beta)\right]\right\} n
\end{equation}

If $\mu(i)=0$ is not taken into account, that is, $s_i=\pm 1$, then according to Kramers and Wannier (1941) or Huang (1987),

\begin{equation}
Q_{_n}=(y+\frac{1}{y})^n=\left[\exp(\beta)+\exp(-\beta)\right]^n
\end{equation}

and,

\begin{equation}\label{eq16}
M(n)\leq\frac{1}{\beta}\log Q_{_n}=\left\{\frac{1}{\beta}\log \left[\exp(\beta)+\exp(-\beta)\right]\right\} n
\end{equation}

On the other hand, if $p$ can be calculated by other methods, we have,

\begin{equation}\label{eq15+1}
M(n)\leq \frac{1}{\beta}\log p Q_{_n}=\frac{1}{\beta}\log p \left[1+\exp(\beta)+\exp(-\beta)\right]^n
\end{equation}

and

\begin{equation}\label{eq16+1}
M(n)\leq\frac{1}{\beta}\log p Q_{_n}=\frac{1}{\beta}\log p\left[\exp(\beta)+\exp(-\beta)\right]^n
\end{equation}

\section{Discussion}

It can be seen from the section \ref{result} that $\mu (n)$ is a particular state of 1D Ising model system in the case of $x=1$, that is, $s_0=s_{_n}=\mu(n)$, $s_1=\mu(1)=+1$, $s_2=\mu(2)=-1$, $s_3=\mu(3)=-1$, $s_4=\mu(4)=0$, $\ldots$, $s_{_{n-1}}=\mu(n-1)$. Its Hamiltonian (total energy) is $\beta \sum \mu(i) = \beta M(n)$, and the probability $p$ for this special can be obtained by eq. (\ref{eq14}). When this Hamiltonian (or $M(n)$) is given or fixed, it can be obtained by numerous ways rather than the sequence of $\mu(n)$ above only, since $s_i$ is an independent random sequence according to the definition of the 1D Ising model. For example, $s_0=s_{_n}=\mu(n)$, $s_1=\mu(10)=+1$, $s_2=\mu(20)=0$, $s_3=\mu(30)=-1$, $s_4=\mu(17)=-1$, $\ldots$, $s_{_{n-1}}=\mu(n-3)$, or, $s_0=s_{_n}=-1$, $s_1=+1$, $s_2=-1$, $s_3=-1$, $s_4=0$, $\ldots$, $s_{_{n-1}}=0$, $\ldots$ ($\sum_{i} s_i=M(n)$), and so on. Therefore, $M(n)$ is equivalent to the sum of an independent random sequence but results in an equal probability $p$ to that from the sequence of $\mu(n)$, i.e.,

\begin{equation}\label{eq19}
p=\frac{\exp\left[\beta M(n)\right]}{Q_{_n}}=\frac{\exp\left[\beta\sum_{i} \mu(i)\right]}{Q_{_n}}=\frac{\exp(\beta\sum_{i} s_i)}{Q_{_n}}
\end{equation}

Thus, 

\begin{equation}\label{eq21}
M(n) =\sum_{i} \mu(i)=\sum_{i} s_i
\end{equation}

Eq. (\ref{eq19}) and (\ref{eq21}) show that $M(n)=\sum_{i} \mu(i)$ is equivalent to $\sum_{i} s_i$   in the sense of equivalent probability in statistical mechanics. 

Therefore, we can estimate the upper bound of the $M(n)$ by the independent random sequence $s_i$ with a expectation $u=0$ and a standard variance $\sigma=\sqrt{2/3}$. Similar to Wei (2016), we have two inequality  
eq. (\ref{eq22}) and (\ref{eq23}) hold with a probability of $1 - \alpha$.

\begin{equation}\label{eq22}
M(n) \le  \sqrt{2/3} K_{_{\frac{\alpha}{2}}}\sqrt{n}
\end{equation}
where,
$$
\int\limits_{ - {K_{\alpha /2}}}^{{K_{\alpha /2}}} {\frac{1}{{\sqrt {2\pi } }}} \exp ( - \frac{{{t^2}}}{2}){\rm{d}}t = 1 - \alpha 
$$
\begin{equation}\label{eq23}
M(n) \le  \frac{\sqrt{2/3}}{\sqrt{\alpha}} \sqrt{n}
\end{equation}

Additionally, we can investigate the upper bound of the $M(n)$ from the energy fluctuations in the canonical ensemble.  From the statistical mechanics we get an Ising system with an energy of $\mathcal{H'}_j(=-\xi h(\sum_{i} s_i)_j)$ has a probability $p_j$,

\begin{equation}\label{eq24}
p_j=\frac{\exp(-\frac{\mathcal{H'}_j}{kT})}{Q_{_n}}=\frac{\exp[\frac{\xi h}{kT}(\sum_{i} s_i)_j]}{Q_{_n}}=\frac{\exp(\beta U_j)}{Q_{_n}}
\end{equation}

The fluctuation of the energy of this system is just the standard deviation of the $U$,

\begin{equation}\label{eq25}
(\Delta U)^2=\sum p_j(U_j-\bar{U})^2=<U^2>-<U>^2
\end{equation}

where $\bar{U}=\sum p_j U_j=<U>$.

On the other hand, 

\begin{equation}\label{eq26}
\frac{\partial\log Q_{_n}}{\partial \beta}=\sum\frac{1}{Q_{_n}} U_j\exp(\beta U_j)=\bar{U}
\end{equation}

and,

\begin{equation}\label{eq27}
\begin{array}{ll}
\frac{\partial^2\log Q_{_n}}{\partial^2 \beta}&=\sum \frac{1}{Q_n} U_j^2\exp(\beta U_j)-\frac{\partial Q_{_n}}{\partial \beta}\sum\frac{1}{Q_{_n}^2} U_j\exp(\beta U_j)\\
&=\sum \frac{1}{Q_n} U_j^2\exp(\beta U_j)-\frac{\partial \log Q_{_n}}{\partial \beta}\sum\frac{1}{Q_{_n}} U_j\exp(\beta U_j)\\
&=<U^2> -<U>^2\\
&=(\Delta U)^2
\end{array}
\end{equation}

Substituting Eq. (\ref{eqQn}) into Eq. (\ref{eq26}) and Eq. (\ref{eq27}) , we get,

\begin{equation}\label{eq28}
\bar{U}=\frac{\partial\log Q_{_n}}{\partial \beta}=\left[\frac{\exp(\beta)-\exp(-\beta)}{1+\exp(\beta)+\exp(-\beta)}\right] n = A\cdot n
\end{equation}

\begin{equation}\label{eq28}
(\Delta U)^2=\frac{\partial^2\log Q_{_n}}{\partial^2 \beta}=\left\{\frac{\exp(\beta)+\exp(-\beta)+4}{[1+\exp(\beta)+\exp(-\beta)]^2}\right\} n = B\cdot n
\end{equation}

Therefore,  let $0<\alpha<1$ the measured value of the $U$ lies with a probability $>1-\alpha$ in the interval 

\begin{equation*}
\left[\bar{U}-\sqrt{\frac{(\Delta U)^2}{\alpha}},\bar{U}+\sqrt{\frac{(\Delta U)^2}{\alpha}}\right]
=\left[A n-\sqrt{\frac{B}{\alpha} n}, A n+\sqrt{\frac{B}{\alpha} n}\right]
\end{equation*}

For the special case of $U_j=\sum \mu(n)=M(n)$, the mean (or expectation) of $\mu(n)$ is 0 (Wei, 2016), and $M(n)$ has a probability $>1-\alpha$ in the interval

\begin{equation*}
\left[-\sqrt{\frac{B}{\alpha} n}, \sqrt{\frac{B}{\alpha} n}\right]
\end{equation*}

That is to say, if $M(n)$ is a measured quantity, its upper bound is $\sqrt{\frac{B}{\alpha} n}$. 

\section{Conclusions}

The upper bound for the Mertens function can be estimated by the way of statistical mechanics, from which the Mertens function can be studied from a physical and statistical perspective. Some new upper bounds, which are shown in the inequality of (\ref{eq15}), (\ref{eq16}), (\ref{eq15+1}) and (\ref{eq16+1}) are presented.

\vspace{5em}


\ \ 

\end{document}